\documentclass[10pt]{article}

 \usepackage{amsfonts,amssymb,graphicx}
\usepackage[colorlinks=true,linkcolor=blue,citecolor=blue]{hyperref}
\usepackage{multicol}

\def\line#1{\hbox to \hsize{#1\hfill}}
\newtheorem{prop}{Proposition}

\newtheorem{defi}{Definition}
\newtheorem{lemm}{Lemma}
\newtheorem{theo}{Theorem}
\newtheorem{coro}{Corollary}
\newtheorem{conj}{Conjecture}

\newcommand{\D}[1][]{\ensuremath{{\mathbb{D}^{#1}} }}
\newcommand{\Z}[1][]{\ensuremath{{\mathbb{Z}^{#1}} }}
\newcommand{\C}[1][]{\ensuremath{{\mathbb{C}^{#1}} }}
\newcommand{\R}[1][]{\ensuremath{{\mathbb{R}^{#1}} }}
\renewcommand{\S}[1][]{\ensuremath{{\mathbb{S}^{#1}} }}
\renewcommand{\H}[1][]{\ensuremath{{\mathbb{H}^{#1}} }}
\newcommand{\T}[1][]{\ensuremath{{\mathbb{T}^{#1}} }}

\newcommand{\G}{\ensuremath{\mathbb{G}}}

\newcommand{\J}{\ensuremath{\mathbb{J}}}

\renewcommand{\L}{{\cal L}}
\newcommand{\N}{{\cal N}}
\newcommand{\M}{{\cal M}}

\newcommand{\co}{{\texttt{cos}\eps}}
\newcommand{\si}{{\texttt{sin}\eps}}
\newcommand{\ex}{{\texttt{exp}\eps}}
\newcommand{\s}{{\cal S}}
\newcommand{\<}{\langle}
\renewcommand{\>}{\rangle}
\newcommand{\ga}{\gamma}
\newcommand{\pa}{\partial}
\newcommand{\de}{\delta}

\newcommand{\eps}{\epsilon}
\newcommand{\te}{\theta}
\newcommand{\ka}{\kappa}
\newcommand{\la}{\lambda}

\newcommand{\sig}{\sigma}

\title{ Hamiltonian stability of Hamiltonian minimal Lagrangian submanifolds  in pseudo- and para-K\"ahler  manifolds}
\author{ Henri Anciaux\footnote{Universidade de S\~ao Paulo; supported by CNPq (PQ 302584/2007-2) and Fapesp (2010/18752-0)}, Nikos Georgiou\footnote{Department of Mathematics and Statistics,
University of Cyprus; partially supported by Fapesp (2010/08669-9)}}

\begin{document}
\maketitle

\centerline{\textbf {\large{Abstract}}}

\bigskip

{\small Let $\L$ be a Lagrangian submanifold of a pseudo- or para-K\"ahler manifold which is H-minimal, i.e.\ a critical point of the volume functional restricted to Hamiltonian variations.
We derive the second variation of the volume of $\L$ with respect to Hamiltonian variations. We apply this formula to several cases. In particular we observe that a minimal Lagrangian submanifold $\L$ in a Ricci-flat pseudo- or para-K\"ahler manifold is H-stable, i.e.\ its second variation
is definite and $\L$ is in particular a local extremizer of the volume with respect to Hamiltonian variations. We also give a stability criterion
for spacelike minimal Lagrangian submanifolds in para-K\"ahler manifolds, similar to Oh's stability criterion for minimal Lagrangian manifolds in K\"ahler-Einstein manifolds (cf \cite{Oh1}). Finally, we determine the H-stability of a series of examples of H-minimal Lagrangian submanifolds: the product $\S^1 (r_1) \times ... \times  \S^1(r_n)$ of $n$ circles of arbitrary radii in complex space $\C^n$ is H-unstable with respect to any indefinite flat Hermitian metric, while the product $ \H^1(r_1) \times ... \times \H^1(r_n)$ of $n$ hyperbolas in para-complex vector space $\D^n$ is H-stable for $n=1,2$ and H-unstable for $n \geq 3.$ Recently, minimal Lagrangian surfaces in the space  of geodesics of space forms  have been characterized (\cite{An3}, see also \cite{Ge}); 
on the other hand, a class of  H-minimal Lagrangian surfaces in the tangent bundle of a Riemannian, oriented surface  has been identified in \cite{AGR}. We discuss the H-stability of all these examples.}

\bigskip

\centerline{\small \em 2000 MSC: 53D12,  49Q05, 53A07
\em }


\section*{Introduction}

 A submanifold  of a pseudo-Riemannian manifold $(\M,g)$ is said to be \em minimal \em if it is a critical point of the volume functional associated to $g.$ The first variation formula  characterizes a  minimal submanifold by the vanishing of the trace of its second fundamental form, the \em mean curvature vector. \em 
 Minimality is therefore the first order condition for a submanifold to be volume-extremizing in its homology class. In \cite{Si}, the second order condition to be satisfied by a volume-extremizing submanifold has been derived in the Riemannian case (see \cite{An2} for the generalization
 in the pseudo-Riemannian setting). Those minimal submanifolds satisfying this condition (hence local extremizers of the volume) are called
 \em stable minimal submanifolds. \em
 
 \medskip
 
 When the ambient manifold $\M$ is in addition symplectic, it is interesting to look at
  those submanifolds on which the symplectic form vanish, namely, \em Lagrangian submanifolds. \em It was discovered in \cite{HL1} that minimal Lagrangian submanifolds
 of a Ricci-flat K\"ahler manifold (usually called \em Calabi-Yau manifolds\em) have the striking property of being \em calibrated. \em A classical argument using Stokes theorem shows that calibrated submanifolds are volume-extremizing, and therefore automatically stable.
 On the other hand,  since the Lagrangian  class is closed
 under the action of Hamiltonian diffeomorphisms, it is natural to pose the problem of extremizing the volume of a Lagrangian submanifold in its Hamiltonian isotopy class. 
A Lagrangian submanifold is said to be \em H-minimal \em (or \em H-stationary\em) if is a critical point of the volume functional with respect
to such Hamiltonian variations. A H-minimal Lagrangian submanifold is characterized by the fact that its mean curvature vector is divergence-free. A H-minimal Lagrangian submanifold will be said to be \em H-stable \em if the second variation of the volume functional, with respect to Hamiltonian variations again, is a definite quadratic form.

\medskip

In \cite{Oh1} and \cite{Oh2}, the second variation formula of a H-minimal Lagrangian submanifold has been derived in the case of a K\"ahler manifold. In the special  case in which $(i)$ the submanifold $\L$ is not only H-minimal  but minimal and $(ii)$ the ambient manifold $\M$ is K\"ahler-Einstein, the second variation operator depends only on the induced metric on $\L$ and it is easy to prove that the $\L$ is H-stable if and only if the first eigenvalue of the Laplacian of the induced metric on $\L$ is larger than the scalar curvature of $\M$ divided by its dimension (\cite{Oh1}). This \em stability criterion \em allowed several classification results: for example,  the Clifford torus is the unique H-stable minimal Lagrangian torus in complex projective plane ${\mathbb C}{\mathbb P}^2$ (\cite{Ha},\cite{Ur1}), and the unique H-stable compact minimal orientable Lagrangian surface of ${\mathbb C}{\mathbb P}^2$ with genus less or equal to $4$ (\cite{Ur2}). There are also some results about minimal Lagrangian submanifolds in the space $L(\S^n)$ of geodesics of the unit sphere $\S^n$ (\cite{Pa}).
On the other hand, it follows from the second variation formula that the product of $r_1 \S^1 \times ... \times r_n \S^1$ of $n$ circles of arbitrary radii in complex Euclidean space $\C^n$, which  is H-minimal but not minimal, is H-stable (\cite{Oh2}). This leads to the so-called \em Oh's conjecture, \em which states that \em  the product of $r_1 \S^1 \times ... \times r_n \S^1$ is a global minimizer of the volume in its Hamiltonian isotopy class. \em In spite of some partial answers (\cite{An1}), the question is still open.

\medskip

  There are two natural extensions of K\"ahler geometry to the pseudo-Riemannian setting: 
a \em pseudo-K\"ahler \em structure $(J,g)$ on a manifold $\M$ is defined by the same axioms than a classical K\"ahler one,  dropping the requirement that the metric $g$ is definite. The fact that $J$ is an isometry implies that the signature of $g$ must be even. Moreover, if the signature of $g$ is $(2p,2(n-p))$, the signature of the induced metric on a Lagrangian submanifold, if non degenerate, must be $(p,n-p).$
On the other hand, a \em para-K\"ahler \em structure is a pair $(J,g)$, with the same properties than a K\"ahler one, except that  $J$ is  para-complex rather than complex, i.e.\ we have $J^2= Id,$ and that the compatibility of $g$ with respect to $J$ becomes 
$g(J.,J.)=-g(.,.)$. The latter equation implies that the signature of $g$ is neutral, but places no restriction on 
the induced metric of a non-degenerate Lagrangian submanifold.  
The simplest example of para-K\"ahler manifold is the Cartesian $n$-product $\D^n$ of the set of para-complex (or \em split-complex, \em or \em double\em) numbers $\D,$
(see Section \ref{SectionOne} for a precise definition of pseudo- and para-K\"ahler structures, as well as  a description of $\D^n$). 

\medskip

Recently there have been growing interest about submanifold theory in the pseudo-Riemannian setting. In \cite{Me} (cf also \cite{HL2}), the concept of calibration is extended to the realm of pseudo-Riemannian geometry. In \cite{Dg} (cf also \cite{An3}), the study of minimal Lagrangian submanifolds
in complex space $\C^n$ endowed with a flat pseudo-Hermitian form is addressed. It also has been observed that the spaces of geodesics of space forms enjoy natural pseudo- or para-K\"ahler structures (cf \cite{GG1}, \cite{AGK}, \cite{An4}).

\medskip

The purpose of this paper is the study of the Hamiltonian stability of H-minimal Lagrangian submanifolds in pseudo- and para-K\"ahler manifolds. In particular, we calculate the second variation of volume with respect to Hamiltonian variations: 

\medskip

\noindent {\bf Main Theorem}: \em Let $(\M,J,g)$ be a pseudo- or para-K\"ahler $2n$-dimensional manifold with Ricci curvature tensor $Ric.$ Let $\L$ be a H-minimal Lagrangian submanifold of $\M$ with mean curvature vector $\vec{H}$ and second fundamental form $h$. Denote respectively by $\Delta$ and $\nabla$ the Laplacian and the gradient of the induced metric on $\L.$
Then, given  $u \in C_c^{\infty}(\L)$ a smooth, compactly supported real function  on $\L$ and $X = J \nabla u$ the corresponding Hamiltonian vector field,
the second variation formula $\delta^2 {\cal V}(\L)(X)$ of the volume of $\L$ with respect to $X$ is given by:
$$\delta^2 {\cal V}(\L)(X)=\int_{\L} \eps \Big((\Delta u)^2 -  Ric(\nabla u,\nabla u)-2 g(n\vec{H},h(\nabla u,\nabla u))\Big)
+g(n\vec{H},J\nabla u)^2, $$
where  $\eps=1$ in the K\"ahler case and $\eps=-1$ in the para-K\"ahler one
 (here and in the following, all the integrals are meant with respect to the volume element
induced by the metric). \em

\medskip

 We observe that our notation differs from that of \cite{Oh2}, since we describe Hamiltonian variations via test functions, instead of using exact $1$-forms as in \cite{Oh2}. 
Of course both viewpoints are equivalent and so in the K\"ahler case we recover Oh's main result.
A first, immediate consequence of this formula is obtained when the last three terms of this formula vanish:

\begin{coro} \label{coroRicciflat} 

Let $(\M,J,g)$ be a Ricci-flat pseudo-or para-K\"ahler manifold and let $\L$ be a minimal Lagrangian submanifold of $\M.$ Then $\delta^2 {\cal V}(\L)(X)=\eps \int_{\L} (\Delta u)^2,$ which implies: 
\begin{itemize}
	\item[--] in the pseudo-K\"ahler case,  $\L$
is a local minimizer of the volume in its Hamiltonian isotopy class;
\item[--]  in the para-K\"ahler case,  $\L$ is a local maximizer of the volume in its Hamiltonian isotopy class.
\end{itemize}
\end{coro}
It has been proved in \cite{Me} (see also \cite{HL2}) in the para-K\"ahler case, and in \cite{An3} in the pseudo-K\"ahler case, that minimal Lagrangian submanifolds of a Ricci-flat manifold enjoy a kind of "Lagrangian calibration" and are therefore extremizers of the volume in their Lagrangian isotopy class. Hence Corollary \ref{coroRicciflat} may be regarded as the local version of the results of \cite{Me} and \cite{An3}. 

\medskip

We also get a result similar to Oh's {\it Stability criterion} for compact, minimal Lagrangian submanifolds with definite induced metric in a para-K\"ahler Einstein manifold: 

\begin{theo}\label{Coro1}
Let $(\M,J,g)$ be a para-K\"ahler Einstein manifold of dimension $2n$ with scalar curvature $2n c$ and $\L$ be a compact,  minimal, Lagrangian submanifold of $\M$
whose induced metric  is definite. Then $\L$ is  a local maximizer of the volume in its Hamiltonian isotopy class if and only if the first 
eigenvalue of the Laplacian of the induced metric on $\L$ satisfies 
$\lambda_1 \geq c.$
\end{theo}

Next, we turn our attention to two natural examples of H-minimal Lagrangian submanifolds which are not minimal. It turns out that the homogeneous tori of complex space $\C^n$, though still H-minimal when the ambient metric is indefinite, are no longer H-stable. Hence there is no counterpart of Oh's conjecture in the indefinite case:

\begin{theo} \label{T^nunstable} Let  $\mathbb T_{r_1,...,r_n}^n:=\S^1(r_1) \times ... \times  \S^1(r_n)$ be the product of $n$ circles of arbitrary radii $r_1, ..., r_n$ in complex space $\C^n$ equipped with the pseudo-K\"ahler structure given by the Hermitian form
$$\<\<.,.\>\>_p := -\sum_{j=1}^{p}  dz_j d\bar{z}_j+\sum_{j=p+1}^{n}  dz_j d\bar{z}_j.$$
Then $\mathbb T_{r_1,...,r_n}^n$ is Lagrangian and H-minimal. Moreover, if $p \neq 0,n$, it is H-unstable. 
\end{theo}

It seems reasonable to expect that Hamiltonian instability is shared by all non-minimal Lagrangian submanifolds, at least in the flat case:

\begin{conj} \label{ConjUnstable}
 Let $\L$ be a H-minimal, non minimal, Lagrangian submanifold in complex space $\C^n$ equipped with the pseudo-K\"ahler structure given by the Hermitian form
$\<\<.,.\>\>_p,$
with $p \neq 0,n.$ Then $\L$ is H-unstable.
\end{conj}

The situation of the para-K\"ahler equivalent of the tori $\T_{r_1,...,r_n}^n$, namely the product of hyperbolas in $\mathbb D^n$, is different and somewhat surprising since here the H-stability depend on their dimension, and not on their induced metric:

\begin{theo} \label{H^nstable} Denote by $\H^{1}_+(r):=\{(x,y) \in \D | \, x^2-y^2=r^2, x > 0 \}$ and 
$\H^1_-(r):=\{(x,y) \in \D | \, x^2-y^2=-r^2, y>0  \}$ the hyperbolas of  "radius" $r$ of the set of para-complex numbers. 
Let  $\mathbb H_{r_1,...,r_n}^n:= \H^1_{\pm}(r_1) \times ... \times \H^1_{\pm}(r_n)$ the product of $n$  hyperbolas of arbitrary signature and radii $r_1, ..., r_n$  in para-complex space $\D^n$ equipped with its canonical para-K\"ahler structure (see Section \ref{SectionOne}). 
Then $\H_{r_1,...,r_n}^n$ is Lagrangian and H-minimal. Morever, the curves $\H^1_{\pm}(r)$ and the surfaces $\H_{r_1,r_2}^2$ are   $H$-stable, while $\H_{r_1,...,r_n}^n, \, n \geq 3,$ is H-unstable.  
\end{theo}

Finally, we compute the stability of some Lagrangian surfaces in two examples of non-flat pseudo- or para-K\"ahler manifolds. First, it has been  recently discovered (\cite{Sa},\cite{GG1},\cite{AGK}) that the spaces of  $L^+(\S^{3}_p)$ and $L^-(\S^{3}_p)$ positive and negative geodesics of $3$-dimensional space form $\S^3_p$ enjoy two natural 
pseudo- or para-K\"ahler structure $(\J,\G)$ and $(\J',\G')$ which share the same symplectic form $\omega=\eps \G(\J.,.)=\eps' \G'(\J'.,.)$. It has been proved in \cite{An4} that $\G$ is Einstein and that $\G'$ is scalar flat, and that Lagrangian surfaces of $L^{\pm}(\S^3_p)$ arise as normal congruences (or \em Gauss maps\em) of surfaces of $\S^3_p$. Moreover, the normal congruence of a tube over a geodesic of $\S^3_p$ is minimal with respect to both metrics $\G$ and $\G'$. The next result discusses the H-stability of such surfaces:

\begin{theo}
Let $\bar{\cal S} \subset L^{\pm}(\S^{3}_p)$ be  the normal congruence of a the tube over a geodesic $\ga$ of $\S^3_p$. Then:  \begin{itemize}
\item[-] $\bar{\cal S}$ is H-stable with respect to $\G$ if and only if $p=2,$ i.e.\ $\S^3_2$ is the anti-de Sitter space and $\ga$ is an unbounded geodesic;
\item[-] $\bar{\cal S}$ is H-stable with respect to $\G'$ if and only if $\ga$ is a closed geodesic (this implies that $p=0$ or $2,$ i.e. $\S^3_p$ is the Riemannian sphere or the anti-de Sitter space).
\end{itemize}
\end{theo}

The last situation we study it the case of the tangent bundle $T\N$  of  a oriented, Riemannian surface $(\N,g_0).$ In  \cite{GK}, it has been proved that  $T\N$  enjoys a natural pseudo-K\"ahler structure $(\J,\G).$ Moreover, it is observed in \cite{AGR} that given a regular curve $\ga$ of $\N,$ the set of tangent vectors to $\N$ which are normal to $\ga$ (thus a surface in $T\N$) is Lagrangian and H-minimal (these surfaces are called \em normal bundles \em in \cite{AGR}). 

\begin{theo}\label{t:theoremontn}
Let  $(\N,g_0)$ be a Riemannian surface with Gaussian curvature $K$,  $\ga$ a  regular curve of $\N$ with geodesic curvature $\ka$  and $\L$ its normal bundle. 
\begin{itemize}
\item[-] If   $\ka^2 \leq -2K$ along $\ga$ (this is for example the case 
if $\ga$ is  a geodesic of a surface $\N$ with negative curvature), then $\L$ is  H-stable;
\item[-] If $\ga$ is closed and $ \sup_{ \ga}(\ka^2+2K) \leq \frac{16 \pi^2}{L^2}, $  where $L$ is the length of $\ga,$ then $\L$ is  H-stable;
\item[-]  If $\ga$ is unbounded and  $\ka^2 > -2K$ along $\ga$ (this is for example the case any curve of a surface with  non-negative curvature), then $\L$ is  H-unstable.
\end{itemize}
\end{theo}

\noindent {\bf Acknowledgements.} The authors would like to thank G. Siciliano for his helpful comments about the proof of Theorems \ref{H^nstable} and \ref{t:theoremontn}.

\section{Preliminaries} \label{SectionOne}

\subsection{The space $\D^n$} \label{Dn}
The set of \em  para-complex \em (or \em split-complex, \em or \em double\em)  numbers
is the two-dimensional real vector space $\R^2$ endowed with the commutative algebra structure whose product  rule is given by
$$ (x , y) . (x' , y')= (xx'+yy', xy'+x'y).$$
The number $(0,1)$, whose square is $(1,0)$ and not $(-1,0)$, will be denoted by $\tau.$
It is convenient to use the following notation: $(x,y) \simeq z =x + \tau y .$
In particular, one has the same conjugation operator than in $\C$:
$$ \overline{x+ \tau y } = x - \tau y,$$
whose corresponding square norm is
$$ | z |^2 :=  z . \bar{z} = x^2 -y^2.$$
In other words, the metric associated to  $|.|^2$ is the Minkowski metric $dx^2 - dy^2.$

\medskip

We define the "para-Cauchy-Riemann" operators on $\D$ by
$$\frac{\pa }{\pa z}:= \frac{1}{2}\left( \frac{\pa }{\pa x} + \tau \frac{\pa }{\pa y} \right) \quad \quad
\frac{\pa }{\pa \bar{z}}:= \frac{1}{2}\left(\frac{\pa }{\pa x} - \tau \frac{\pa }{\pa y} \right).$$
A map $f$ defined on a domain of $\D$ is said to be \em para-holomorphic \em if it satisfies
$\frac{\pa f}{\pa \bar{z}}=0.$
Observe that this is a hyperbolic equation, so a para-holomorphic map needs not to be analytic, and may be merely continuously differentiable.

\medskip

On the Cartesian $n$-product $\D^n$ with para-complex coordinates $(z_1, ...,z_n)$, we define the canonical  para-K\"ahler structure $(J_*,\<.,.\>_*)$ on $\D^n$ by
$$ J_*(z_1,...,z_n):=(\tau z_1, ...,\tau z_n)$$
$$ \<.,.\>_*:=\sum_{j=1}^{n} dz_j  d\bar{z}_j = \sum_{j=1}^{n}  dx_j^2-dy_j^2.$$
We also introduce the "para-Hermitian" form:
$$\<\<.,.\>\>_* := \sum_{j=1}^{n} dz_j \otimes d\bar{z}_j= \sum_{j=1}^{n}  dx_j^2-dy_j^2 - \tau  \sum_{j=1}^{n} dx_j \wedge dy_j.$$
In other words we recover $\<.,.\>_*$ by taking the real part of  $\<\<.,.\>\>_*$, while its imaginary part is minus the associated symplectic form
$\omega:=\<J_*.,.\>_*.$

\subsection{Pseudo- and para-K\"ahler manifolds} 
Let $\M$ be a $2n$-dimensional smooth manifold.  An almost complex (resp.\ almost para-complex) structure $J$ is a $(1,1)$-tensor
such that $J^2=- Id$ (resp.\ $J^2 = Id$). In the almost para-complex
case we furthermore require that the two eigen-distributions $Ker(J \pm Id)$ have the same rank\footnote{In the almost complex case, the two
eigen-distributions $Ker(J \pm i Id)$ of $T\M \otimes \C$ have always the same rank.}. Since we shall always deal simultaneously with the two cases, 
we set $\eps$ equal to $1$ in the complex case and $-1$ in the para-complex case, so that we have $J^2=-\eps Id.$
 Next we introduce the \em Nijenhuis tensor:\em
$$N^J(X,Y):=[X,Y]-\eps([JX,JY]-J[JX,Y]-J[X,JY]).$$
The tensor $N^J$ vanishes if and only if $J$ is actually a complex (resp.\ para-complex structure), i.e.\ there exists an atlas on $\M$ whose transition maps are local bi-holomorphic (resp.\ bi-para-holomorphic) diffeomorphisms of $\C^n$ (resp.\ $\D^n$)\footnote{Observe that in the complex case, this is due to quite a hard result of Newlander and Nirenberg (see \cite{NN}), while in the para-complex case it is a simple consequence of Frobenius theorem.}.
A pseudo-Riemannian metric $g$ is said to be \em compatible \em with $J$ if 
$$g(J.,J.) = \eps g.$$

A pair $(J,g)$, where $J$ is pseudo- or para-complex structure, and $g$ is a compatible pseudo-Riemannian metric, is said to be \em a K\"ahler structure \em on $\M$ if the alternated $2$-from $\omega:=\eps g(J.,.)$ is symplectic, i.e.\ $d\omega=0.$ One may roughly summarize this by saying that pseudo- or para-K\"ahler geometry is the intersection of complex (or para-complex), pseudo-Riemannian and symplectic geometries.

\subsection{Volume minimization with respect to Hamiltonian variations}

\begin{defi}
Let $(\M,\omega)$ a $2n$-dimensional symplectic manifold. 

\begin{itemize}

\item[--]
A submanifold $\L$ of $\M$ is said to be \em Lagrangian \em if it has dimension $n$ and if the symplectic form $\omega$ vanishes on it.  
\item[--]
A vector field $X$ on $\M$ is said to be \em Hamiltonian \em if the $1$-form $X \lrcorner \,  \omega$ is exact.
\end{itemize}

\end{defi}
The relevance of Hamiltonian vector fields in the study of Lagrangian submanifolds is explained in the following proposition:
\begin{prop}
Let $X$ be a Hamiltonian vector field on $\M.$ Then the Lie derivative of $\omega$ along $X$ vanishes, and therefore if $\L$ is Lagrangian, so does $\exp_t(X)(\L).$ 
\end{prop}

\begin{prop}
$(\M,J,g)$ be a pseudo- or para-K\"ahler manifold, $\L$ a non-degenerate Lagrangian submanifold  and $X$ a compactly supported Hamiltonian vector field whose restriction to
$\L$ is normal to $\L.$ Then there exists $u \in C^{\infty}_c(\L)$ such that $X|_{\L}= J \nabla u,$ where $\nabla$ denotes the gradient operator on $\L$ with respect to the induced metric.
\end{prop}

 We recall that the first variation formula states that
$$
\delta {\cal V}(\L)(X):=\left. \frac{d}{dt}{\cal V}(\exp_t(X)(\L))\right|_{t=0}=-\int_{\L} g(n\vec{H},X)dv.
$$
(see the next section for the definition of the mean curvature $\vec{H}$).

\begin{coro}
Let $(\M,J,g)$ be a pseudo- or para-K\"ahler manifold. Then a Lagrangian non-degenerate submanifold $\L$ is H-minimal, i.e.\ a critical point of the volume restricted to Hamiltonian variations, if and only if its mean curvature vector $\vec{H}$ satisfies the equation 
$$ div (nJ\vec{H})= 0,$$
where $div$ denotes the divergence operator on $\L$ with respect to the induced metric. 
\end{coro}

\begin{defi}
Let $(\M,J,g)$ be a pseudo- or para-K\"ahler manifold. A H-minimal Lagrangian submanifold
 $\L$ is said to be \em H-stable  \em if  the quadratic form restricted to Hamiltonian  variations
$\delta^2 {\cal V}(\L)(X):= \left. \frac{d^2}{dt^2}{\cal V}(\exp_t(X)(\L))\right|_{t=0}$ is definite.
\end{defi}

\subsection{The second fundamental form of a Lagrangian submanifold}
We recall briefly the definition of the second fundamental form of a submanifold $\L$ of $(\M, J,g)$ with non degenerate induced metric.
At a point $x \in \L \subset \M,$ we have
 the splitting $T_x \M = T_x \L \oplus N_x \L $ and for $X \in T_x \M ,$ we write $ X = X^\top + X^\bot$ accordingly.
We recall that the second fundamental form of $\L$ is defined by
$h(X,Y)=(\nabla_X Y)^\bot, $  where ${\nabla}$ denotes the Levi-Civita connection of $g$, and its shape operator by $A_N X = -(\nabla_X N )^\top,$ where $\nabla^\bot$
is the connection of the normal bundle $N \L.$  These two tensors carry actually the same information (the \em extrinsic geometry \em of $\L$) since they are related by the formula
$ g(h(X,Y),N) = g(A_N X, Y).$ The \em mean curvature \em $\vec{H}$ of $\L$ is the normal vector field obtained by tracing the second fundamental form $h$ with respect to the induced metric $g.$

\medskip

If $\L$ is in addition Lagrangian, the formula $\omega=\eps g(J.,.)$ implies the $J$ is an isometry (or anti-isometry) of the tangent bundle $T\L$ of $\L$ onto its normal bundle $N\L.$ This fact implies  
the following symmetry property:
\begin{lemm} \em (\cite{An2}, \cite{An4}) \em \label{tri} Let $\L$ be a non-degenerate, Lagrangian submanifold of a pseudo-K\"ahler or para-K\"ahler
manifold $(\M,J,g).$ 
Then
the map $$g(h(X, Y),J Z)=g(A_{JZ} X,Y)$$ is
tri-symmetric, i.e.\ $${g}(h(X,Y),JZ)=g(h(Y,X),JZ)={g}(h(X,Z),JY).$$

\end{lemm}

\section{The second variation formula under Hamiltonian variations (Proof of the Main Theorem)}

 We start from the evolution equation satisfied by the volume density $dv$ associated to the induced metric
$$\frac{d}{dt}(dv)=-g(n\vec{H},X)dv,$$
which, in particular, implies the well-known first variation formula
\begin{equation}\label{e:firstvariation111}
\frac{d}{dt}{\cal V}(\exp_t(X)(\L))=-\int_{\exp_t(X)(\L)}g(n\vec{H},X)dv.
\end{equation}
Differentiating  Equation \ref{e:firstvariation111} at $t=0$ yields
\begin{equation} \label{secondvariation}
\delta^2{\cal V}(\L)(X)=-\int_{\L}\frac{d}{dt} \left.\Big(g(n\vec{H},X)\Big)\right|_{t=0}dv+\int_{\L} g(n\vec{H},X)^2dv.
\end{equation}
Letting $(e_1,...,e_n)$ be a frame which is orthonormal at $t=0,$ we set $g_{ij}:=g(e_i,e_j)$ and $g^{ij}$ to be the coefficients of the inverse matrix of $[g_{ij}]_{1 \leq i,j \leq n}.$
Using the Einstein sum convention, we have
\begin{eqnarray}
\frac{d}{dt} g(n\vec{H},X)&=& \frac{d}{dt} \left( g^{ij}g(h(e_i,e_j),X) \right) \nonumber\\
&=&\frac{dg^{ij}}{dt} g(h(e_i,e_j),X)+g^{ij}g \left(\frac{d}{dt}h(e_i,e_j),X \right)+g^{ij}g \left(h(e_i,e_j),\frac{dX}{dt} \right) \nonumber  \\
&=&\frac{dg^{ij}}{dt} g(h(e_i,e_j),X)+g^{ij}g \left(\frac{d}{dt}h(e_i,e_j),X \right)+g \left( n\vec{H},\frac{dX}{dt} \right) \label{e:oitreisoroi}.
\end{eqnarray}
We deal with the first two terms of Equation \ref{e:oitreisoroi} in the following Lemma:

\begin{lemm}\label{l:firsttwotermsofscvr}
We have:
$$  \left. \frac{dg^{ij}}{dt}\right|_{t=0} g(h(e_i,e_j),X)
+g^{ij}g \left( \left. \frac{d}{dt}h(e_i,e_j)\right|_{t=0},X \right)=-g(\nabla^{\bot}X,\nabla^{\bot}X)-g(R^{\bot}X,X)+g(A_{X},A_{X})+ \Delta f,$$
where $f:=\frac{1}{2}g(X,X).$
\end{lemm}

\noindent
{\it Proof.} The proof is elementary and can be found  in \cite{An2} (Proof of Theorem 5, Chapter 1)\footnote{The Reader is warned that although it is assumed in \cite{An2} that $\vec{H}$ vanishes, this assumption
is not used in the computation of the formula, which is therefore still valid here.} or in \cite{Xi}).

\medskip

The  third term of Equation \ref{e:oitreisoroi} does not vanish since $\L$ is not  minimal, nor the vector field $\frac{dX}{dt}$ is Hamiltonian  \em a priori\em:
\begin{lemm}\label{l:vectorhamiltonian}
We have:
\[
\int_{\L}   g \left(n\vec{H},\;\frac{dX}{dt}\Big|_{t=0} \right) =\eps \int_{\L} g(n\vec{H},h(\nabla u,\nabla u)) .
\]
\end{lemm}

\noindent
{\it Proof.} Since the vector field $X(.,t)$ is Hamiltonian  $\forall t,$ there exists $t_0 >0$  and $\bar{u}(x,t)\in C^{\infty}_{c}(\L \times (-t_0,t_0))$ such that $\bar{u}(.,0)=u$ and 
$X(.,t)=J\nabla \bar{u}(.,t).$
If follows that
\begin{eqnarray*}
 \frac{dX}{dt}\Big|_{t=0} &=&  \frac{d}{dt}g^{ij} (e_i(\bar{u}) Je_j)  \Big|_{t=0} \\
&=&\frac{dg^{ij}}{dt}\Big|_{t=0}e_i(u)Je_j+\epsilon_i\frac{de_i(\bar{u})}{dt}\Big|_{t=0}Je_i  +\eps_i e_i(u)  \frac{d}{dt} J e_i |_{t=0} .
\end{eqnarray*}
Taking into account the fact that $\frac{dg^{ij}}{dt}\Big|_{t=0}=2\eps_{i}\eps_j g(h(e_i,e_j),X)$ 
(cf Lemma 4, p.\ 24, \cite{An2})
and the tri-symmetry of the curvature tensor $ g(h(.,.),J.)$ (Lemma \ref{tri}), we first calculate
\begin{eqnarray*}\frac{dg^{ij}}{dt}\Big|_{t=0}e_i(u)Je_j &=&    2\eps_{i}\eps_j g(h(e_i,e_j), J \nabla u) e_i(u)Je_j  \\
&=&  2\eps_{i}\eps_j \eps_k g(h(e_i,e_j), J e_k)  e_k(u) e_i(u)Je_j \\
&=&  2\eps_{i}\eps_j \eps_k g(h(e_i,e_k), J e_j)  e_k(u) e_i(u)Je_j \\
&=& 2 \eps_j g(h(\nabla u, \nabla u), Je_j) Je_j\\
&=& 2 \eps h(\nabla u, \nabla u).
\end{eqnarray*}
On the other hand,   using the fact that $[X,e_i]=0$ at $t=0$, we get
\begin{eqnarray*}
\eps_i e_i(u)  \frac{d}{dt} J e_i \Big|_{t=0} &=&\eps_i e_i(u) \nabla_X J e_i \\
&=&   \eps_i e_i(u) J \nabla_{e_i}^{\bot} X\\
&=& - \eps_i e_i(u) J A_X e_i\\
&=& - J A_X \nabla u,
\end{eqnarray*}
so that, using again Lemma \ref{tri}:
\begin{eqnarray*}   g \left( n\vec{H} , \eps_i e_i(u)  \frac{d}{dt} J e_i  \Big|_{t=0} \right) &=&-   g \left( n\vec{H} , J A_{J \nabla u} \nabla u \right)\\
&=&  g( J n\vec{H},   A_{J\nabla u} \nabla u)\\
&=&    g(  h(J n\vec{H}, \nabla u ), J \nabla u) \\
&=&  g(  h( \nabla u, \nabla u ), J (J n \vec{H}) ) \\
&=& -\eps  g(  h( \nabla u, \nabla u ), n \vec{H}) ).
\end{eqnarray*}
We conclude, taking into account that the vector field $ J \nabla \left(\left. \frac{\pa \bar{u} }{\pa t} \right|_{t=0} \right)$ is Hamiltonian,
\begin{eqnarray*}
\int_{\L} g \left( n\vec{H},\frac{dX}{dt} \Big|_{t=0} \right)&=&
\int_{\L} g\left(n\vec{H},\frac{dg^{ij}}{dt}\Big|_{t=0}e_i(u)Je_j+\epsilon_i\frac{de_i(\bar{u})}{dt}\Big|_{t=0}Je_i  +\eps_i e_i(u)  \frac{d}{dt} J e_i |_{t=0} \right)\\
&=& \int_{\L}  2\eps g \left( n\vec{H}  , h(\nabla u, \nabla u) \right) + g \left(  n\vec{H}  ,J \nabla \left( \frac{\pa \bar{u} }{\pa t} \Big|_{t=0} \right) \right) 
                                                            -\eps  g(  h( \nabla u, \nabla u ), n \vec{H}) )\\
&=&  \eps   \int_{\L}   g(n\vec{H},h(\nabla u, \nabla u)), 
\end{eqnarray*}
 which is the claimed formula.


\bigskip

Using Lemmas \ref{l:firsttwotermsofscvr}, \ref{l:vectorhamiltonian} and Equation \ref{e:oitreisoroi}, the second variation formula  (Equation \ref{secondvariation}) becomes  
\begin{equation}\label{e:constaruction12}
\delta^2{\cal V}(\L)=\int_{\L} g(\nabla^{\bot}X ,\nabla^{\bot}X) +g(R^{\bot}X,X)-g(A_X,A_X)-\eps g(h(\nabla u,\nabla u),n\vec{H})+g(n\vec{H},X)^2.
\end{equation}
The next Lemma relates two terms of the above expression to the Ricci tensors of $\M$ and $\L$:
\begin{lemm}\label{l:importlemma}
Let $V$ be a tangent vector to $\L$, then the following formula holds:
$$g(R^{\bot}JV,JV)-g(A_{JV},A_{JV})= \eps \Big(-Ric^{\M}(V,V)+ Ric^{\L}(V,V)- g(n\vec{H},h(V,V))  \Big)$$
where $Ric^{\M}$ denotes the Ricci tensor of $(\M,g)$ and $Ric^{\L}$ denotes the Ricci tensor of $\L$ endowed with the induced metric. 
\end{lemm}

\noindent {\it Proof:}
Observe first that
\begin{eqnarray*}
  g( A_{JV},A_{JV})&=&   \sum_{i,j=1}^n  \eps_i \eps_j  g( A_{JV}e_i, e_j)\\
&=&  \sum_{i,j=1}^n  \eps_i \eps_j  g(h(e_i,e_j), JV)\\
&=&  \sum_{i,j=1}^n  \eps_i \eps_j  g(h(V,e_i), Je_j)  \\
&=& \eps  \sum_{i,j=1}^n  \eps_i (\eps \eps_j)  g(h(V,e_i), Je_j)  \\
&=&  \eps \sum_{i=1}^n   g(h(V,e_i),h(V,e_i)) .
\end{eqnarray*}
We conclude, using Gauss equation for immersions and the fact that
$Ric^{\M}(J.,J.)=\eps Ric^{\M}(.,.),$
\begin{eqnarray*}
g(R^\perp JV , JV  )&=& \sum_{i=1}^n \eps_i g(R^\M (JV,e_i)JV, e_i) \\
&=& Ric^M (JV,JV) - \sum_{i=1}^n \eps \eps_i g(R^\M(JV,Je_i)JV, Je_i)\\
&=& \eps Ric^\M (V,V) -  \sum_{i=1}^n \eps \eps_i g(R^\M(V,e_i)V, e_i)\\
&=& \eps Ric^\M (V,V)  - \sum_{i=1}^n \eps \eps_i \left( g( R^\L(V,e_i)V, e_i) - g(h(e_i,e_i), h(V,V))+ g(h(V,e_i),h(V,e_i))  \right) \\
&=&   \eps Ric^\M (V,V)- \eps Ric^\L (V,V)  - \eps g( n\vec{H}, h(V,V)) +    g( A_{JV},A_{JV}).
\end{eqnarray*}


\bigskip

We observe now that $g(\nabla^{\bot} J\nabla u, \nabla^{\bot} J\nabla u)  =\eps g(\nabla^2 u,\nabla^2 u),$ where the Hessian operator $\nabla^2 u$ is defined by
$ \nabla^2 u(X,Y)=g(\nabla_X \nabla u, Y ) = g(\nabla_Y \nabla u, X ).$ Using Lemma \ref{l:importlemma}, Equation \ref{e:constaruction12}, we get
\begin{equation}\label{e:secvarform}
\delta^2{\cal V}(\L)(J \nabla u)=\int_{\L} \eps \Big( g(\nabla^2 u,\nabla^2 u)-Ric^{\M}(\nabla u,\nabla u)+Ric^{\L}(\nabla u,\nabla u)-2 g( n\vec{H}, h(\nabla u,\nabla u))\Big)+ g(n\vec{H},J\nabla u)^2.
\end{equation}

The next step consists of extending a formula due to Reilly (see \cite{Re}) to the pseudo-Riemannian setting:
\begin{lemm}\label{l:nicolas}
Let $(\L,g)$ be a pseudo-Riemannian manifold with Ricci tensor $Ric^\L$ and $u$ a smooth, compactly supported function on $\L.$ Then
$$ \int_{\L} (\Delta u)^2 - g(\nabla^2 u,\nabla^2 u)  = \int_{\L} Ric^\L(\nabla u,\nabla u).$$
\end{lemm}

\noindent
{\it Proof.} The proof is based on the generalization of Bochner's formula to the pseudo-Riemannian setting whose proof is postponed to the Appendix of this paper:
$$\frac{1}{2}\Delta (g(\nabla u,\nabla u))=Ric^\L(\nabla u,\nabla u)+g(\nabla u,\nabla(\Delta u))+g(\nabla^2 u,\nabla^2 u).$$

Integrating this equation over $\L$ yields
\begin{eqnarray*}
0 \, =\int_{\L}\frac{1}{2}\Delta(g(\nabla u,\nabla u))&=&\int_{\L}Ric^\L(\nabla u,\nabla u)+ g(\nabla u,\nabla(\Delta u))+g(\nabla u^2,\nabla u^2) \nonumber \\
&=&\int_{\L}Ric^\L(\nabla u,\nabla u)-(\Delta u)^2+g(\nabla u^2,\nabla u^2).
\end{eqnarray*}


We are now in position to complete the proof of the Main theorem: applying Lemma \ref{l:nicolas} to Equation \ref{e:secvarform} gives
$$\delta^2 {\cal V}(\L)(X)=\int_{\L} \eps \Big((\Delta u)^2 -  Ric^\L(\nabla u,\nabla u)-2 g(n\vec{H},h(\nabla u,\nabla u))\Big)
+g(n\vec{H},J\nabla u)^2 ,$$
the required formula.


\section{Applications}

\subsection{Minimal Lagrangian submanifolds in para-K\"ahler-Einstein manifolds (Proof of Theorem \ref{Coro1})}
In this section we assume that the Lagrangian submanifold $\L$ is not only H-minimal but also minimal and that $(\M,g)$ is Einstein, i.e.\ there exists a real function $c$ such that
$Ric^{\M} = c g.$ It is  well known that $c$ must be constant and  that the scalar curvature of $(\M,g)$ is $ 2n c.$
 Then the second variation formula  of the volume of $\L$ with respect Hamiltonian variations becomes:
 \begin{eqnarray*}
\delta^2{\cal V}(\L)(J\nabla u)&=&\eps \int_{\L} (\Delta u)^2 - Ric^{\M}(\nabla u,\nabla u) \\
 &=&\eps \int_{\L} (\Delta u)^2 - c g(\nabla u,\nabla u)\\
 &=&\eps \int_{\L}  (\Delta u)^2 + c u  \Delta u .
\end{eqnarray*}
Observe that this formula depends only on the induced metric on $\L.$ 

\medskip

If we assume furthermore that  $\L$ is compact and that its induced metric is definite, we may use the  spectral properties of the Laplacian  $-\Delta$ of a Riemannian metric: 
it is well known that  $-\Delta$ is an elliptic operator and that the space $L^2(\L)$ of square integrable real functions of $\L$ enjoys a Hilbertian basis  $(\phi_i)_{i \in \mathbb N}$ 
satisfying $-\Delta \phi_i = \la_i \phi_i,$ where $(\la_i)_{ i \in \mathbb N}$ is an increasing sequence of strictly positive numbers tending to infinity.
Hence, given $u = \sum_{i \in \mathbb N} a_i \phi_i \in C^\infty_c(\L) \subset L^2(\L)$, we have
 $-\Delta u =  \sum_{i \in \mathbb N} \la_i a_i \phi_i.$
If follows that
\begin{eqnarray*}
\delta^2{\cal V}(\L)(J \nabla u)&=&\eps \big( \<\Delta u, \Delta u\>_{L^2(\L)} + c \<u, \Delta u\>_{L^2(\L)}  \big)\\
&=&\eps \left( \sum_{i \in \mathbb N} \lambda_i^2 a_i^2  - c \sum_{i \in \mathbb N} \lambda_i a_i^2 \right)\\
&=& \eps  \sum_{i \in \mathbb N} a_i^2 \lambda_i (\lambda_i - c).
\end{eqnarray*}
Therefore, if $ \lambda_1 \geq c,$ the second variation $\delta^2{\cal V}(\L)(J \nabla u)$ has the sign of $\eps$ and $\L$ is H-stable.

\vspace{0.1in}

Conversely,  if $\lambda_1 < c,$ it easy to see that
$$\eps  \delta^2{\cal V}(\L)( J \nabla \phi_1) = \lambda_1(\lambda_1 - c) <0,$$
 while, since $\lim_{i \to \infty} \lambda_i = \infty, $ there exists $i_0$ such that
$$\eps \delta^2{\cal V}(\L)( J \nabla \phi_{i_0}) = \lambda_{i_0}( \lambda_{i_0} - c) >0.$$
 Hence $\L$ is H-unstable.

In the para-K\"ahler case, we get Theorem \ref{Coro1}. On the other hand, if  $(\M,J,g)$ is pseudo-K\"ahler with signature $(2n,2(n-p))$, the induced metric on $\L$ has signature $(n,n-p).$ Hence the discussion above applies only in the K\"ahler case. We therefore recover exactly Oh's stability criterion:

\medskip

\noindent \textbf{Stability Criterion (\cite{Oh1})}: \em Let $(\M, J,g)$ be an Einstein-K\"ahler manifold with
with scalar curvature $2nc$ and $\L$ is a compact, minimal, Lagrangian submanifold of $\M.$ Let $\lambda_1$ be
the first eigenvalue of the Laplacian. Then $\L$ is H-stable
 if and only if $\lambda_1>c.$ \em

\subsection{Tori $\mathbb T_{r_1,...,r_n}^n$ are H-unstable (proof of Theorem \ref{T^nunstable})} \label{secunstable}
We endow the complex vector space $\C^n$ with the Hermitian form
$$\<\<.,.\>\>_p :=\sum_{j=1}^n \eps_j dz_j \otimes d\bar{z}_j= \<.,.\>_p - i \omega_p,$$
where $\eps_j = \pm 1.$
 A parametrization of the torus $\T_{r_1,...,r_n}^n:=\S^1(r_1) \times ... \times \S^1(r_n) \subset \C^n$ is
 $$ 
  \begin{array}{lccc}  f :
  & \R / 2\pi r_1 \Z\times ... \times \R / 2\pi r_n \Z
  &\to&  \C^n \\
 & (s_1,...,s_n)& \mapsto &  (r_1 \exp(is_1/r_1), ... , r_n \exp(i s_n/r_n)).
 \end{array}
$$
The first derivatives of $f$ are
$$f_{s_j}  =  i \exp(is_j/r_j) e_j , $$
where $(e_1,...,e_n)$ is the canonical Hermitian basis of $\C^n.$
It is straightforward to check that $f$ is Lagrangian with respect to $\omega_p$ and that the induced metric $f^*\<.,.\>_{2p}$
takes the form $\sum_{j=1}^n \eps_j ds_j^2$. In particular it is non-degenerate and flat.

Moreover, setting
$$ N_j := J f_{s_j} = -  \exp(is_j/r_j)e_j,  $$
we have
$$h_{jkl}:=\<f_{s_j s_k},N_l\>_{2p}=\frac{\eps_j \de_{jk} \de_{kl}}{r_j} .$$

Therefore 
$$\<n\vec{H}, N_j\>_{2p} = \sum _{k=1}^n \eps_k h_{jkk}=\frac{1}{r_j}$$
and
$$ n \vec{H}=  \sum _{j=1}^n \frac{\eps_j}{r_j}N_j .$$
In particular, ${div} (nJ \vec{H})$ vanishes, i.e.\ $f$ is H-minimal.

On the other hand, 
letting $u$ be a smooth map on $\T_{r_1,...,r_n}^n,$
we have $\nabla  u =  \sum_{j=1}^n \eps_j u_{s_j} \pa_{s_j} ,$
so $J \nabla u=- \sum_{j=1}^n \eps_j u_{s_j} N_j.$
It follows that
$$\<n\vec{H},J \nabla u \>_{2p}= -\sum_{j=1}^n\frac{u_{s_j}}{r_j}$$
and
\begin{eqnarray*} \<h(\nabla u, \nabla u ),n \vec{H} \>_{2p}&=& \sum_{j=1}^n u_{s_j}^2 \<h(\pa_{s_j}, \pa_{s_j} ),2 \vec{H} \>_{2p}\\
 &=&  \sum_{j=1}^n u_{s_j}^2 \frac{h_{jjj}}{r_j}\\
 &=& \sum_{j=1}^n u_{s_j}^2 \frac{u_{s_j}^2}{r_j^2}. 
\end{eqnarray*}
Finally, applying the Main Theorem,
\begin{eqnarray*} \delta^2{\cal V}(\T_{r_1,...,r_n}^n)(J\nabla u)&=&
\int_{\R / 2\pi r_1 \Z \times ... \times \R / 2\pi r_n \Z}   \left( \left(\sum_{j=1}^n \eps_j u_{s_j s_j}\right)^2  - 2 \sum_{j=1}^n \frac{u_{s_j}^2}{r_j^2}+  
  \left( \sum_{j=1}^n \eps_j \frac{u_{s_j}}{r_j}  \right)^2  \right)ds_1 ... ds_n \\
     & =& \int_{\R / 2\pi r_1 \Z \times ... \times \R / 2\pi r_n \Z} \left( \left(\sum_{j=1}^n \eps_j u_{s_j s_j}\right)^2  - \sum_{j=1}^n \frac{u_{s_j}^2}{r_j^2}   +  2\sum_{j <k} \eps_j \eps_k \frac{u_{s_j} u_{s_k}}{r_j r_k}   \right)ds_1 ... ds_n. \end{eqnarray*}
       
From this expression, it is easy to check that the second variation is indefinite, exhibiting both positive and negative directions: 

\begin{itemize}
\item[$\bullet$]On the one hand, taking $u(s_1,...,s_n)=\cos (k  s_1/r_1 )$, where $k$ is an arbitrary integer, we get
\begin{eqnarray*}
\delta{\cal V}^2(\T_{r_1,...,r_2}^n)(J\nabla u) &=&
 ( \prod_{j=2}^n 2\pi r_j) \left(\frac{k^4 }{r_1^4}\int_0^{2\pi r_1} (\cos (k s_1/r_1))^2 ds_1 - \frac{k^2 }{r_1^4}\int_0^{2\pi r_1}(\sin (k s_1/r_1))^2 ds_1\right)\\
 &=& (\prod_{j=2}^n 2\pi r_j) \frac{\pi(k^4-k^2)}{r_1^3}
\end{eqnarray*}
which is strictly positive if  $k >1$.

\item[$\bullet$]
On the other hand, since
$p \neq 0, n$, there exists both positive and negative vectors in the canonical basis $(e_1, ...,e_n).$ Assume, without loss of generality, that $\eps_1=1$ and $\eps_2=-1.$
Consider first the case  $\frac{r_1}{r_2} \in \mathbb Q$; we shall use the fact that the wave equation 
$u_{s_1 s_1}-u_{s_2 s_2}=0$ admits globally defined solutions: given a $2\pi$-periodic, non constant
map $F$, and two integers $a$ and $b$ satisfying $\frac{a}{b}=\frac{r_1}{r_2}$, we set
$ u(s_1,...,s_n):= F(\frac{as_1}{r_1}-\frac{b s_2}{r_2}).$
We then get
$$ \sum_{j=1}^n \eps_j u_{s_j s_j}= u_{s_1 s_1}-u_{s_2 s_2}=\left( \frac{a^2}{r_1^2}-\frac{b^2}{r_2^2} \right) F''=0$$
and
$$-\sum_{j=1}^n \frac{u_{s_j}^2}{r_j^2}+  2\sum_{j <k} \eps_j \eps_k \frac{u_{s_j} u_{s_k}}{r_j r_k} 
=-\left( \frac{u_{s_1}}{r_1} + \frac{u_{s_2}}{r_2}\right)^2
= -\left(\frac{a}{r_1^2}+\frac{b}{r_2^2} \right)^2 (F')^2.$$
Hence
$$ \delta^2{\cal V}(\T_{r_1,...,r_2}^n)(J\nabla u) =- \prod_{j=3}^n  (2\pi r_j) \left(\frac{a}{r_1^2}+\frac{b}{r_2^2} \right)^2 
\int_{\R / 2\pi r_1 \Z \times \R / 2\pi r_2 } (F')^2 ds_1 ds_2 <0.$$

To conclude, observe  we have proved H-instability of $\T^n_{r_1,...,r_n}$  for a dense set of values of $\frac{r_1}{r_2}.$ 
The H-instability  being an open condition, it follows that all tori $\T^n_{r_1,...,r_n}$ are H-unstable.
\end{itemize}

\subsection{Stability and instability of the products  $\H_{r_1,...,r_n}^n$ (proof of Theorem~\ref{H^nstable})} \label{secunstable}
We refer to Section \ref{Dn} for notation and introduce  furthermore:
$$
 \ex(\tau t):= \left\{ 
\begin{array}{ccc} \cosh (t) + \tau \sinh (t) &\mbox{ if }& \eps=1, \\
\sinh (t) + \tau \cosh (t) &\mbox{ if } &\eps=-1.\end{array} \right. 
$$
Observe that $\ex(\tau t)'=\texttt{exp}(-\eps)(\tau t)=\tau \ex(\tau t)$
and that $|\ex(\tau t)|^2=\eps.$
 A parametrization  of $\H_{r_1,...,r_n}^n$ is:
 $$ \begin{array}{lccc}  f :
  &  \R^n 
  &\to&  \D^n \\
 & (s_1,...,s_n)& \mapsto &  (r_1 \ex_1(\tau s_1/r_1), ...,  r_n \ex_n(\tau s_n/r_n)).
 \end{array}$$
The first derivatives of $f$ are
$$f_{s_j}  =  \tau \ex_j ( \tau s_j / r_j)e_j , $$
where $(e_1,...,e_n)$ is the canonical para-Hermitian basis of $\D^n.$
It is straightforward to check that $f$ is Lagrangian with respect to $\omega_*$
and that the induced metric $f^*\<.,.\>_*$ takes the form
 $-\sum_{j=1}^n \eps_j ds_j^2$. In particular it is non-degenerate and flat.

Next, setting
$ N_j := J_* f_{s_j} =  \ex_j ( \tau s_j / r_j)e_j  ,$
we easily calculate the mean curvature vector of $f$:
$$ n \vec{H}= -\sum_{j=1}^n \frac{\eps_j}{r_j}N_j .$$
In particular, ${div} (nJ_* \vec{H})$ vanishes, i.e.\ $f$ is H-minimal.

On the other hand, 
letting $u$ be a compactly supported smooth map on $\H_{r_1,...,r_n},$
we have $\nabla  u = - \sum_{j=1}^n \eps_j u_{s_j} \pa_{s_j} ,$
so $J_* \nabla u= -\sum_{j=1}^n \eps_j u_{s_j} N_j.$
It follows that
$$\<n\vec{H},J_* \nabla u \>_*= \sum_{j=1}^n\frac{\eps_j u_{s_j}}{r_j}$$
and 
\begin{eqnarray*} \<h(\nabla u, \nabla u ),n \vec{H} \>_*&=& \sum_{j=1}^n u_{s_j}^2 \<h(\pa_{s_j}, \pa_{s_j} ),n \vec{H} \>_*\\
 &=& -\sum_{j=1}^n u_{s_j}^2 \frac{\eps_j h_{jjj}}{r_j} \\
 &=& -\sum_{j=1}^n \frac{u_{s_j}^2}{r_j^2}. \end{eqnarray*}
\bigskip

Hence, by the Main Theorem,
\begin{eqnarray*} \delta^2{\cal V}(\H^n_{r_1,...,r_n})(J_* \nabla u)&=&
\int_{\R^n}\left(- \left( \sum_{j=1}^n\eps_j u_{s_j s_j} \right)^2  -2 \sum_{j=1}^n \frac{u_{s_j}^2}{r_j^2}+   \left(\sum_{j=1}^n \eps_j \frac{u_{s_j}}{r_j} \right)^2  \right)ds_1 ... ds_n\\
     & =& \int_{\R^n} \left( -\left( \sum_{j=1}^n \eps_j u_{s_j s_j} \right)^2  - \sum_{j=1}^n \frac{u_{s_j}^2}{r_j^2} +  2\sum_{j <k} \eps_j \eps_k \frac{u_{s_j} u_{s_k}}{r_j r_k}  \right)ds_1 ... ds_n .
      \end{eqnarray*}
   In the one-dimensional case, this expressions becomes
  $$ \delta^2{\cal L}(\H^{1}_{\pm}(r_1))(J_* \nabla u)=
 \int_{\R} \left( - (u'')^2   - \left(\frac{u'}{r_1} \right)^2 \right)ds_1  \leq 0,$$
 while in the two-dimensional case, we have
$$\delta^2{\cal A}(\H^2_{r_1,r_2})(J_* \nabla u)=
 \int_{\R^2} \left( -\left( \eps_1 u_{s_1 s_1} + \eps_2 u_{s_2 s_2} \right)^2  - \left(\eps_1 \frac{u_{s_1}}{r_1} -
  \eps_2 \frac{u_{s_2} }{r_2} \right)^2 \right)ds_1 ds_2  \leq  0.$$
 Hence in dimension $1$ and $2$, the products $\H^n_{r_1,...,r_n}$ are H-stable. We now prove that it is not anymore the case
if $n \geq 3.$ This follows from
three elementary lemmas.
For sake of brevity we set $s=(s_1,...,s_n)$ and $ds=ds_1 ... ds_n.$ 

\begin{lemm} \label{firstlemmaQ}
Set
$$ Q(\nabla u , \nabla u):=  \sum_{j=1}^n \frac{u_{s_j}^2}{r_j^2} -  2\sum_{j <k} \eps_j \eps_k \frac{u_{s_j} u_{s_k}}{r_j r_k}.$$
If the quadratic functional
$${\textbf Q}(u):= \int_{\R^n} Q(\nabla u , \nabla u) ds$$
is indefinite, then the second variation 
$$  \delta^2{\cal V}(\H^n_{r_1,...,r_n})(J_* \nabla u)=
 -\int_{\R^n}  \left( \Delta u \right)^2 ds -  \int_{\R^n}  Q(\nabla u , \nabla u)   ds$$
 is indefinite 
as well.
\end{lemm}

\noindent
{\it Proof.} Given $u \in C^{\infty}_c(\R^n)$, set $u^t(s):=t^{n/2-1} u(t s).$
A quick calculation, using that $ u^t_{s_i}=t^{n/2}  u_{s_i}$ and $u^t_{s_i s_i}=t^{n/2+1}u_{s_i s_i},$ 
shows that
$$  \delta^2{\cal V}(\H^n_{r_1,...,r_n})(J_* \nabla u^t)=  - t^{2} \int_{\R^n}  \left( \Delta u \right)^2   ds +     \int_{\R^n} Q(\nabla u , \nabla u)   ds.$$
Letting $t$ tends to $0$ makes this expression to have the same sign than ${\textbf Q}(u),$
 so the indefiniteness of $\textbf{Q}$ implies that of $ \delta^2{\cal V}(\H^n_{r_1,...,r_n}).$

\begin{lemm}
There exists a linear change of variable $\sig(s)$ and an integer $k$, $1 \leq k \leq n-1,$ such that 
$${\textbf Q} (u):=c \int_{\R^n}  \left(\sum_{i=1}^k u_{\sigma_i}^2 - \sum_{i=k+1}^n u_{\sigma_i}^2 \right) d\sig,$$
where $c$ is a non vanishing constant.
\end{lemm}

{\it Proof.}  Observe that
$$Q(\nabla u, \nabla u)= [ \nabla u ]. M .   [\nabla u ]^{\rm{T}} ,$$
where the matrix $M$ is defined by 
$$ M:= \mbox{diag}(1/r_1^2, ..., 1/r_n^2) -\left[\frac{\eps_i \eps_j }{r_i r_j} \right]_{1 \leq i,j \leq n}.$$ 
On the one hand, we have
$$ (1,0,...,0) . M . ( 1,0,...,0)^{\rm{T}} = \frac{1}{r_i^2} > 0$$
and on the other hand, $-(n-2)$ is an eigenvalue of $M$, with eigenvector $V =(\frac{\eps_1}{r_1}, ... , \frac{\eps_n}{r_n}),$ so
$$ V . M . V^{\rm{T}}= -(n-2) |V|_0^2<0.$$
Therefore the matrix $M$ is indefinite. The conclusion follows from Sylvester's law of inertia.

\begin{lemm}
The quadratic functional
$${\textbf Q} (u):= \int_{\R^n}  \left(\sum_{i=1}^k u_{\sig_i}^2 - \sum_{i=k+1}^n u_{\sig_i}^2 \right) d\sig$$
 is indefinite on the space of compactly supported, smooth functions $C^{\infty}_c(\R^n).$
\end{lemm}

{\it Proof.} 
 Given a non vanishing $u \in C^{\infty}_c(\R^n)$,  set $u^t(\sig_1,...,\sig_n):= u(t \sig_1, ...,t \sig_k ,  \sig_{k+1}, ...,  \sig_n)$.
A quick calculation
shows that
$${\textbf Q} (u^t):=t^{-k} \Big(t^2 \int_{\R^n}  \big(\sum_{i=1}^k u_{\sigma_i}^2 \big) d\sig  - \int_{\R^n}  \big( \sum_{i=k+1}^n u_{\sigma_i}^2 \big) d\sig \Big).$$
Since both integrals are non zero, letting $t$ tend to $\infty$ and $0$ give respectively a positive and negative value to ${\textbf Q}.$ 
  
\subsection{Minimal Lagrangian surfaces in the space of oriented geodesics of 3-dimensional space forms}
We briefly recall the construction of the canonical pseudo- or para-K\"ahler structures of the space of geodesics of the 3-dimensional space forms.
For further detail, see \cite{An4}.
Consider first the flat pseudo-Riemannian metric $\<\cdot ,\cdot\>_p$ of $\R^4$ of signature $(p,4-p)$:
\[
\<\cdot ,\cdot\>_p:=-\sum_{i=1}^pdx_i^2+\sum_{i=p+1}^4dx_i^2,
\]
and the $3$-dimensional quadric
\[
\S^3_p=\{x\in \R^4 \big|\<x,x\>_p=1\}.
\]

 The space $L^+(\S^3_p)$ (resp. $L^{-}(\S^3_p)$)  of positive (resp.\ negative) oriented geodesics of the space form $\S^3_p$ can be identified with the Grassmannian $Gr^{+}(4,2)$ (resp.\  $Gr^{+}(4,2)$) of oriented two-planes of $\R^4$ with positive induced metric (resp.\ indefinite induced metric). Such Grassmannians, which are 4-dimensional,   may be naturally embedded in  the space of bi-vectors, a $6$-dimensional real linear space: 
\[
 \iota:L^{\pm}(\S^3_p)\rightarrow \Lambda^2(\R^4)=Span\{e_i\wedge e_j\; :\; 1\leq i<j\leq 4\},
\]
where $(e_1,e_2,e_3,e_4)$ denotes the canonical basis of $\R^4$. We endow $\Lambda^2(\R^4)$  with the flat pseudo-Riemannian metric
\[
\<\<x\wedge y, x'\wedge y'\>\>=\<x,x'\>_p\<y,y'\>_p-\<x,y'\>_p\<y,x'\>_p.
\]
A tangent vector to $\iota(L^{\pm}(\S^3_p))$ at $x\wedge y$ takes the form $x\wedge X+y\wedge Y$, where $X,Y$ belong to the orthogonal $ (x\wedge y)^{\bot}$  of  $x\wedge y.$ 
Any of the two  2-dimensional, oriented planes  $x\wedge y$ and $ (x \wedge y)^\perp$ is endowed with a natural complex or para-complex structure that we shall denote by ${\rm J}$ and ${\rm J}'$ respectively. For example ${\rm J}$ is defined by  ${\rm J}x:=y$ and ${\rm J}y:=-\eps x,$ where  $\eps := \<y,y\>_p.$ 

\smallskip

We are now in position to define the two structures $(\J,\G)$ and $(\J',\G').$ We first set $\G$ to be the induced metric on $L^{\pm}(\S^3_p)$, i.e., ${\mathbb G}:=\iota^{\ast}\<\<.,. \>\>$. 
Next we set ${\mathbb J}$ and $\J'$ to be respectively:
$${\mathbb J}(x\wedge X+y\wedge Y):={\rm J} x\wedge X+{\rm J}y\wedge Y.$$
and
$${\mathbb J}'(x\wedge X+y\wedge Y):=x\wedge {\rm J}' X+y\wedge {\rm J}' Y.$$
Finally, the neutral metric $\G'$ is defined by
$$\G':=\epsilon \G(\J.,\J'.):=-\eps  \G(.,\J \circ \J'.) .$$
This implies that the two pseudo- or para-K\"ahler structures $(\J,\G)$ and $(\J',\G')$ share the same symplectic form $\omega=\eps \G(\J.,.)=\eps' \G'(\J'.,.).$ Moreover  the normal congruence of an oriented surface $\s$ of $\S_p^3,$ i.e.\ the set of geodesics normal to $\s,$
is a Lagrangian surface, and if $\s$ is an equidistant tube over a geodesic, its normal congruence is minimal with respect to both the metric $\G$ and $\G'$ (see \cite{An4}). 
In the remainder of the section we shall discuss the H-stability of such surfaces with respect to $\G$ and $\G'.$

\medskip

Let $\ga$ be a non-null geodesic of $\S^3_p$, with arclength parameter $s$ and $(n_1,n_2)$ a local, orthonormal frame of its normal bundle (hence $(\ga,\ga',n_1,n_2) \in SO(p,n-p)$) such that  $D_{\ga'}n_i = n_i(s)'=0.$ The topological type of the tube over $\ga$ depends on the type of the geodesic $\ga$ and on the causal character of its normal space. In order to simplify the exposition, we set $\eps_1:=\<\ga', \ga'\>_p$
and $\eps_4:=\<n_1,n_1\>_p \<n_2,n_2\>_p$ (it will be clear in a moment that this apparently  unnatural notation is chosen to match that of \cite{An4}).
Then the tube over $\ga$ is topologically  $\S^1_{\eps_1} \times \S^1_{\eps_4},$  
where
$$\S^1_{\eps}:=\{ (\co (\sig),\si(\sig) ) | \, \sig \in \R\},$$
i.e.  a circle if $\eps=1$ and a hyperbola if $\eps=-1.$
Hence a local parametrization of the tube over $\ga$ of radius $\te$ is
  is
 $$ \begin{array}{lccc}  \phi :
  & \S^1_{\eps_1} \times \S^1_{\eps_4}  
  &\to&  \S^3_p \\
 & (s,t)& \mapsto &  \co(\te) \ga + \si (\te)\big(\co_4(t) n_1 + \si_4(t) n_2 \big).
 \end{array}$$
Introducing 
$$e_1 :=\ga' \quad \quad e_2:=-\eps_4 \si_4(t) n_1 + \co_4(t) n_2,$$
  we have
\begin{eqnarray*} \phi_s&=& \co(\te) \ga'= \co(\te) e_1\\
\phi_t&=& \si (\te)\big(-\eps_4 \si_4(t) n_1 + \co_4(t) n_2 \big)= \si (\te) e_2.
\end{eqnarray*}
On the other hand,
a unit normal vector is given by
$$N= -\eps \si(\te) \ga +  \co (\te)\big(\co_4(t) n_1 + \si_4(t) n_2 \big).$$
Observe that $\<N,N\>_p=\<n_1,n_1\>_p= \eps.$
Next we have
\begin{eqnarray*} N_s&=& -\eps \si (\te) \ga' =  -\eps \si (\te)  e_1\\
N_t&=&   \co (\te)\big(-\eps_4 \si_4(t) n_1 + \co_4(t) n_2 \big) = \co (\te) e_2 .
\end{eqnarray*}
The normal congruence of $\phi$ is parametrized by $\bar{\phi}:=\phi \wedge N$ (see \cite{An4}).
Setting
\begin{eqnarray*}E_1 := \phi \wedge e_1 &\quad \quad& E_2 := \phi \wedge e_2\\
E_3 := N \wedge e_1 &\quad \quad& E_4 := N \wedge e_2,
\end{eqnarray*}
we have
$$\bar{\phi}_s=\phi_s \wedge N + \phi \wedge N_s=-\co (\te) E_3 - \eps \si(\te) E_1$$
$$\bar{\phi}_t=\phi_t \wedge N + \phi \wedge N_t=-\si (\te) E_4+ \co (\te) E_2.$$
Since the coefficients of the metrics $\G$ and $\G'$ in the basis $(E_1, E_2, E_3, E_4)$
 are
$$ \G = \mbox{diag}  ( \eps_1, \eps' \eps_1 , \eps \eps_1 , \eps \eps' \eps_1):= ( \eps_1, \eps_2,  \eps_3 , \eps_4)$$
and
$$ \G' =\left( \begin{array}{cccc} 0 & 0 & 0 & \eps_2 \\ 0 & 0 & -\eps_2 & 0 \\ 0 & -\eps_2 &0&0 \\ \eps_2 & 0 & 0 &0 \end{array} \right),$$
the coefficient of induced metrics $\bar{\phi}^*(\G)$ and $\bar{\phi}^*(\G')$ in the coordinates $(s,t)$ are
$$ \left( \begin{array}{cc}
\eps_3 & 0 \\ 
0 & \eps_2
\end{array} \right) \quad \mbox{ and }  \quad \left( \begin{array}{cc}
 0&  \eps_2 \\ \eps_2 &0 
\end{array} \right). $$
Moreover, the coefficients of the restriction of the Ricci tensor  $\bar{\phi}^*(Ric^\G)$ and $\bar{\phi}^*(Ric^{\G'})$ in the coordinates $(s,t)$ are
$$ 2\left( \begin{array}{cc}
\eps_1 & 0 \\ 
0 & \eps_4
\end{array} \right) \quad \mbox{ and }  \quad -2 \left( \begin{array}{cc}
 \eps_1& 0 \\ 
0 & \eps_4
\end{array} \right) .$$

We deduce that the Hamiltonian second variation of the area of $\bar{\s}$ with respect to $\G$ is
\begin{eqnarray*}
\delta^2{\cal A}_{\G}(\bar{\s})(\J \nabla u)&=&\eps \int_{\bar{\s}} \left(({\Delta}_{g} u)^2-Ric^{\G}(\nabla u,\nabla u) \right) \\
&=&\eps \int_{\S^1_{\eps_1} \times \S^1_{\eps_4}}\left( (\eps_3 u_{ss}+ \eps_2 u_{tt})^2-2(\eps_1 u^2_s+\eps _4 u^2_t) \right) ds dt,
\end{eqnarray*}
while the Hamiltonian second variation of the area of $\bar{\s}$ with respect to $\G'$ is
\begin{eqnarray*}
\delta^2 {\cal A}_{\G'}(\bar{\s})(\J' \nabla u)&=&\eps' \int_{\bar{\s}} \left((\Delta_{g'} u)^2-Ric^{\G'}(\nabla u,\nabla u) \right) \\
&=&\eps' \int_{\S^1_{\eps_1} \times \S^1_{\eps_4}} \left(4u^2_{st}+2(\eps_1 u^2_s+ \eps_4 u^2_t) \right) ds dt .
\end{eqnarray*}

\subsubsection{The case of $L(\S^3)$}
In this case we have $(\eps_1,\eps_2,\eps_3,\eps_4)=(1,1,1,1)$, so
$$ \delta^2{\cal A}_{\G}(\bar{\s})(\J \nabla u) =\int_{\S^1 \times \S^1}\left( (u_{ss}+ u_{tt})^2- 2(u^2_s+u^2_t) \right) ds dt.$$
Since $(L^+(\S^3),\J,\G)$ is K\"ahler-Einstein with scalar curvature $8$ and that the first eigenvalue of the Laplacian on $\S^1 \times \S^1$ endowed with the flat metric is $1,$
we deduce that that $\bar{\s}$ is H-unstable with respect to $\G$. One could also make a direct proof, 
getting both a positive and a negative second variation considering $u_1=\cos(s+t)$ and $u_2 = \cos (2s)$, so the second variation is
 indefinite\footnote{The H-instability of $\bar{\s}$ with respect to $\G$ was proved in \cite{Pa}.}. On the other hand,
$$ \delta^2{\cal A}_{\G'}( \bar{\s})(\J' \nabla u) =\int_{\S^1 \times \S^1} \left( 4u_{st}^2+ 2(u^2_s+u^2_t) \right) ds dt,$$
which  is positive.

\subsubsection{The case of an unbounded geodesic of $L(d\S^3)$} \label{unboudeddS3}
Here the metric is $\<.,.\>_1:= -dx_1^2 + dx_2^2+dx_3^2+dx_4^2.$ Hence an unbounded geodesic is, modulo congruence,
$\ga(s)=(\sinh s, \cosh s ,0,0),$
so $\ga'=(\cosh s, \sinh s,0,0)$ and $\eps_1=-1.$
Hence we may choose for the normal vectors $n_1=(0,0,1,0)$ and $n_2=(0,0,0,1)$, which implies $\eps=1$ and $\eps_2=1.$ Hence $\eps'=-1$ and $\eps_4=1.$
Finally $(\eps_1,\eps_2,\eps_3,\eps_4)=(-1,1,-1,1)$, so we get
$$ \delta^2{\cal A}_\G( \bar{\s})(\J \nabla u) =\int_{\R \times \S^1}\left( (-u_{ss}+  u_{tt})^2+ 2u^2_s-2u^2_t \right) ds dt.$$
We claim that the latter is indefinite: on the one hand, we have $ \delta^2{\cal A}_{\G}( \bar{\s})(\J \nabla a(s))=2\pi \int_{\R} ((a'')^2+ 2(a')^2)ds >0$ and on the other hand
  $\delta^2{\cal A}_{\G}( \bar{\s})(\J \nabla (a(s)\cos t))=\pi \int_{\R} ((a'')^2+ 4(a')^2 - a^2)ds $
which may take negative values.

\medskip
Analogously, we obtain
$$ \delta^2{\cal A}_{\G'}( \bar{\s})(\J' \nabla u) =-\int_{\R \times \S^1}\left( 4u_{st}^2- 2u^2_s+u^2_t \right) ds dt,$$
which is indefinite as well, since, for example $ \delta^2{\cal A}_{\G'}( \bar{\s})(\J' \nabla a(s))=4\pi \int_{\R} (a')^2ds >0$ and \\
$ \delta^2{\cal A}_{\G'}( \bar{\s})(\J' \nabla (a(s)\cos t))=- \pi \int_{\R} (2(a')^2+  a^2)ds<0. $

\subsubsection{The case of a definite tube over a closed geodesic of $L(d\S^3)$}
A closed geodesic is, modulo congruence,
$\ga(s)=(0,0,\cos s, \sin s),$ 
so $\ga'=(0,0,-\sin s, \cos s)$ and $\eps_1=1.$
In order to get a tube with definite metric, we choose $n_1=(1,0,0,0)$ and $n_2=(0,1,0,0)$, which yields $\eps=-1$ and $ \eps_2=1.$ Hence $\eps'=1$ and $\eps_4=-1.$
Finally $(\eps_1,\eps_2,\eps_3,\eps_4)=(1,1,-1,-1)$, so we get
$$ \delta^2{\cal A}_\G( \bar{\s})(\J \nabla u) =-\int_{\S^1 \times \R}\left( (-u_{ss}+  u_{tt})^2- 2u^2_s+2u^2_t \right) ds dt,$$
and
$$ \delta^2{\cal A}_{\G'}( \bar{\s})(\J' \nabla u) =\int_{\S^1 \times \R}\left( 4u_{st}^2+ 2u^2_s- 2u^2_t \right) ds dt,$$
which are both indefinite, as in the previous subsection.

\subsubsection{The case of an indefinite tube over a closed geodesic of $L(d\S^3)$}
Again we have
$\ga(s)=(0,0,\cos s, \sin s),$ so $\ga'=(0,0,-\sin s, \cos s)$ and $\eps_1=1.$
In order to get a tube with indefinite metric, we choose $n_1=(0,1,0,0)$ and $n_2=(1,0,0,0)$, which yields $\eps=1$ and  $\eps_2=-1.$ Hence $\eps'=-1$ and $\eps_4=-1.$
It follows that
$$ \delta^2{\cal A}_\G( \bar{\s})(\J \nabla u) =\int_{\S^1 \times \R}\left( (u_{ss}-  u_{tt})^2- 2u^2_s+2u^2_t \right) ds dt,$$
and
$$ \delta^2{\cal A'}_{\G'}( \bar{\s})(\J' \nabla u) =-\int_{\S^1 \times \R}\left( 4u_{st}^2+ 2u^2_s-2u^2_t \right) ds dt.$$
As above, both quadratic forms are indefinite.

\subsubsection{The case of a closed  geodesic of $L(Ad\S^3)$}
In the same way as above, we obtain
$$ \delta^2{\cal A}_{\G}( \bar{\s})(\J \nabla u) =-\int_{\S^1 \times \S^1}\left( (u_{ss}+  u_{tt})^2- 2u^2_s-2u^2_t \right) ds dt,$$
which is indefinite (exactly as in the spherical case), and
$$ \delta^2{\cal A}_{\G'}( \bar{\s})(\J' \nabla u) =-\int_{\S^1 \times \S^1}\left( 4u_{st}^2+ 2(u^2_s+u^2_t) \right) ds dt,$$
 which is negative.

\subsubsection{The case of an indefinite tube over an unbounded closed  geodesic of $L(Ad\S^3)$}


Proceeding as in the previous sections, we obtain
$$ \delta^2{\cal A}_{\G}( \bar{\s})(\J \nabla u) =-\int_{\R^2}\left( (u_{ss}+  u_{tt})^2+ 2(u^2_s+u^2_t) \right) ds dt,$$
which is clearly negative, and
$$ \delta^2{\cal A}_{\G'}( \bar{\s})(\J' \nabla u) =-\int_{\R^2}\left( 4u_{st}^2- 2(u^2_s+u^2_t) \right) ds dt,$$
 which is indefinite, by an argument similar to that of Lemma \ref{firstlemmaQ} of Section \ref{secunstable}: setting  $u^\lambda(s,t):= u(\lambda s, \lambda t),$ we have
$$ \delta^2{\cal A}_{\G'}( \bar{\s})(\J' \nabla u^\lambda) =- 4 \lambda^2 \int_{\R^2} u_{st}^2ds dt + 2   \int_{\R^2}(u^2_s+u^2_t)  ds dt.$$

\subsubsection{The case of a definite tube over an unbounded closed  geodesic of $L(Ad\S^3)$}


Here we have
$$ \delta^2{\cal A}_{\G}( \bar{\s})(\J \nabla u) =\int_{\R^2}\left( (u_{ss}+  u_{tt})^2+ 2(u^2_s+u^2_t) \right) ds dt,$$
which is positive,
$$ \delta^2{\cal A}_{\G'}( \bar{\s})(\J' \nabla u) =\int_{\R^2}\left( 4u_{st}^2- 2(u^2_s+u^2_t) \right) ds dt,$$
 which is indefinite as it has been seen in the previous subsection.

\subsubsection{The case of $L(\H^3)$}

We have
$$ \delta^2{\cal A}_{\G}( \bar{\s})(\J \nabla u) =-\int_{\R \times \S^1}\left( (u_{ss}-  u_{tt})^2+ 2u^2_s-2u^2_t \right) ds dt,$$
which is indefinite (see Section \ref{unboudeddS3}.)
and
$$ \delta^2{\cal A}_{\G'}( \bar{\s})(\J' \nabla u) =\int_{\R \times \S^1}\left( 4u_{st}^2- 2u^2_s+2u^2_t) \right) ds dt,$$
which is indefinite as well.

\bigskip

\bigskip

\bigskip

 \hspace{-7em}  \begin{tabular}{| l || l || c| c|  c | c | c |}
     \hline
     Space form &  Type of $\gamma$ & Induced metric on $\s$
     &  $(\eps_1,\eps_2, \eps_3,\eps_4)$ & Topology  of $\s$ & $\G$-H-stability    &  $\G'$-H-stability \\ \hline \hline
   $\S^3_0=\S^3 $ & closed & definite   & $(1,1,1,1)$& torus  & unstable & stable  \\ \hline 
   $ \S^3_1=d\S^3$ &  closed  & definite &  $ (1,1,-1,-1)$ & cylinder & unstable & unstable \\ \hline 
 & closed  & indefinite &  $ (1,-1,1,-1) $ & cylinder & unstable & unstable \\ \hline 
   & unbounded & indefinite &  $(-1,1,-1,1)$ & cylinder  & unstable & unstable \\  \hline 
  $\S^3_2=Ad\S^3$ & closed & indefinite  &  $(1,-1,-1,1)$ & torus  &  unstable & stable \\ \hline 
   &  unbounded & indefinite & $(-1,1,1,-1)$ & plane  & stable & unstable \\ \hline 
&  unbounded & definite   & $ (-1,-1,-1,-1)$ & plane  & stable & unstable \\ \hline 
  $\S^3_3=\H^3$   &  unbounded & definite &  $(-1,-1,1,1)$ &   cylinder &  unstable & unstable  \\ 
     \hline
   \end{tabular}

\subsection{Rank 1 Lagrangian surfaces in the tangent bundle of a Riemann surface}\label{s:tangentbundlofn}

Let  $(\N ,g_0)$ be   an oriented Riemannian surface with metric $g_0$ and denote by $j$ the  complex structure induced from its orientation (i.e.\ the rotation of angle $\pi/2$). It has been proved in \cite{GK} (cf also \cite{AGR},\cite{An2}) that the tangent bundle of $\N$ enjoys a natural pseudo-K\"ahler structure $(\J,\G)$, which may be described as follows: the Levi-Civita connection of $g$ induces a splitting of $TT\N$ into 
\begin{eqnarray*}
TT\N=H\N\oplus V\N&\simeq& T\N\oplus T\N  \\
X &\simeq & ( X^h,X^v),
\end{eqnarray*}
According to this decomposition, we set
\begin{eqnarray*}
	 \J (X^h,X^v)&:=&(j X^h, j X^v) \\
	 \G\big((X^h,X^v),(Y^h,Y^v)\big)&:=& g_0(  X^v, jY^h)-g_0( X^h,jY^v).
	 \end{eqnarray*}
It has been proved in \cite{GK} (see also \cite{AR}) that  $(T\N,\J,\G)$ is a pseudo-K\"ahler manifold and that $\G$ have neutral signature. The next proposition describes a class of H-minimal Lagrangian surfaces of $(T\N,\J,\G)$:

\begin{prop}[\cite{AGR}]

Let $\ga: I \to \N$ be a regular curve parametrized by arclength and a real map $a \in C^1(I).$ Then the surface $ \L$ image of the immersion
\begin{eqnarray}
f: I \times \R &\rightarrow& T\N\nonumber \\
(s,t)&\mapsto & (\gamma(s),\; a(s)\gamma'(s)+tj\gamma'(s)),\nonumber
\end{eqnarray}
is a H-minimal, flat, Lagrangian surface.  Moreover, it is minimal if and only if  $\gamma$ is a geodesic. 
These surfaces are characterized by the fact that the restriction to $\L$ of the canonical projection $\pi: T\N \to \N$ has rank one.
\end{prop}
In order  to study the H-stability of such surfaces (Theorem \ref{t:theoremontn}), we shall need the following description of 
the Ricci curvature of the metric $\G$:

\begin{prop}[\cite{AR}] \label{curvature} 
Let $X$ and $Y$ and  two vector fields on $T\N$ and $K$ the Gaussian curvature of $(\N,g_0)$. Then the Ricci curvature of $(T\N,\G)$
is given by
$${Ric}^{\G}(X,Y)=  2 Ric^{g_0} (X^h, Y^h) = 2K g_0(X^h,  Y^h).$$
\end{prop}

\medskip

\noindent {\it Proof of the Theorem \ref{t:theoremontn}}: 
An easy calculation (cf \cite{AGR}) shows that the coefficients of the
induced metric $f^{\ast}\G$ in the coordinates $(s,t)$ are
$$ \left( \begin{array}{cc}
-2a\kappa & -1 \\ 
-1 & 0
\end{array} \right).$$
If follows that the gradient of the induced metric takes the following form
$$\nabla u = -u_t \pa_s +(-u_s+2 a \ka u_t)\pa_t.$$
On the other hand, the coefficients of the second fundamental form $h$ are
$$h_{112}= \ka \quad \quad h_{122}=0 \quad \mbox{ and } \quad h_{222}=0,$$
so in particular the mean curvature vector is $2 \vec{H}=\ka \J \pa_t.$ 
Hence
$$
\G(2\vec{H},\J \nabla u)=-\G(\J 2 \vec{H},\nabla u )
=- \G(-\ka \pa_t, -u_t \pa_s +(-u_s+2 a \ka u_t)\pa_t)=- \G(\ka \pa_t, u_t \pa_s)= \ka u_{t}.
$$
Next we have
\[
h(\nabla u,\nabla u)=u_{t}^2h(\pa_{s},\pa_{s})-2u_{t}(2a\ka u_{t}-u_{s})h(\pa_t,\pa_s)+(2a\ka u_{t}-u_{s})^2h(\pa_t,\pa_t),
\]
so
\begin{eqnarray*}
 \G(h(\nabla u,\nabla u),2 \vec{H} )&=&  \ka  \G(h(\nabla u,\nabla u),\J \pa_t )\\
 &=&  \ka  \left(u_{t}^2 h_{112}-2u_{t}(-u_s+2a\ka u_t)h_{112}+(-u_s+2a\ka u_t)^2 h_{222} \right)\\
 &=&\ka^2 u_{t}^2.
\end{eqnarray*}
 By Proposition \ref{curvature} we have
\[
Ric^{\G}(\pa_s,\pa_s)=2K\qquad\quad Ric^{\G}(\pa_{s},\pa_{t})=Ric^{\G}(\pa_{t},\pa_{t})=0,
\]
so
$$
Ric^{\G}(\nabla u,\nabla u)=2Ku_{t}^2.
$$
 By the Main Theorem, we obtain
$$\delta^2{\cal A}(\L)(\J \nabla u)=\int_{I \times \R}\Big(4u_{st}^2- (\ka^2+2K)u^2_{t}\Big)ds dt.$$
In particular, if  $\ka^2 \leq -2K$ along the curve $\ga$, the Hamiltonian second variation is obviously positive.

\medskip

We now claim that if  $\ga$ is unbounded and $\ka^2 > -2K$, then $\L$ is H-unstable. To see this, 
set $u^\la(s,t):= \la^{3/2} u (\la s,t)$, so that
 $$\delta^2{\cal A}(\L)(\J \nabla u^\la)=\la^2 \int_{\R^2}4u_{st}^2 ds dt  - \int_{\R^2} \big(\ka^2(\la^{-1}s)+2K(\la^{-1}s) \big)u^2_{t}ds dt,$$
so letting $\la$ tend to $0$ and $\infty$ respectively yields positive and negative variations.

\medskip

Finally, we assume that $\ga$ is a closed curve of length $L.$ Then, 
using Wirtinger's inequality for $L$-periodic functions, we have
\begin{eqnarray*}
\delta^2{\cal A}(\L)(\J \nabla u)&=&\int_{\R / L \Z \times \R}\Big(4u_{st}^2-(\ka^2+2K)u^2_{t}\Big)ds dt\\
& \geq & \frac{4 \pi^2}{L^2}  \int_{\R / L \Z \times \R} 4 u_{t}^2 ds dt-  (\sup_{\R / L \Z} (\ka^2+ 2K))  \int_{\R / L \Z \times \R} u^2_{t}ds dt\\
 &\geq &   \left(\frac{16 \pi^2}{L^2}-\sup_{\R / L \Z} (\ka^2+ 2K) \right)  \int_{\R / L \Z \times \R} u^2_{t}ds dt.
\end{eqnarray*}
Hence, if 
$ \sup_\ga(\ka^2+2K) \leq \frac{16 \pi^2}{L^2},  $
then $\L$ is H-stable.

\section*{Appendix}

The purpose of this Appendix is the proof of the following:
\medskip

\noindent \textbf{Pseudo Riemannian Bochner's formula}
\em Let $(\L,g)$ be a pseudo-Riemannian manifold with Ricci tensor $Ric^\L$ and $u$ a smooth, compactly supported function on $\L.$ Then  \em
$$
\frac{1}{2}\Delta (g(\nabla u,\nabla u))=Ric^\L(\nabla u,\nabla u)+g(\nabla u,\nabla(\Delta u))+g(\nabla^2 u,\nabla^2 u).
$$

\noindent
{\it Proof.} 

\medskip

 Let   $(e_1, ... ,e_n)$ be a local orthonormal frame which is normal at some point $x,$ i.e. $\nabla_{e_i}e_j(x)$ vanishes. This implies that
$\nabla_{\nabla u} e_i=\sum \eps_i e_i(u) \nabla_{e_i}e_j(x)$ vanishes as well.
We set as usual $\eps_i:=g(e_i,e_i).$ Then we have
\begin{eqnarray*}
\frac{1}{2}\Delta g(\nabla u,\nabla u)&=&\frac{1}{2}\sum_{i=1}^n\eps_i   e_i (e_i (g(\nabla u,\nabla u)))\\
 &=& \sum_{i=1}^n\eps_i   e_i (g(\nabla_{e_i} \nabla u,\nabla u))\\
&=&   \sum_{i=1}^n\eps_i   e_i ( \nabla^2 u(e_i, \nabla u))\\
&=&   \sum_{i=1}^n\eps_i   e_i (g( \nabla_{\nabla u} \nabla u,e_i)))\\
&=&   \sum_{i=1}^n\eps_i   g ( \nabla_{e_i} \nabla_{\nabla  u}  \nabla u, e_i)) \\
&=&   \sum_{i=1}^n\eps_i   \left(g ( R({\nabla  u},e_i)  \nabla u, e_i)) + g( \nabla_{\nabla  u} \nabla_{e_i}\nabla u, e_i) - g(\nabla_{[\nabla u ,e_i]} \nabla u,e_i) \right)  \\
&=& Ric^\L(\nabla u,\nabla u) +   \sum_{i=1}^n\eps_i   \left( g( \nabla_{\nabla  u} \nabla_{e_i}\nabla u, e_i) - g(\nabla_{[\nabla u ,e_i]} \nabla u,e_i) \right) .
\end{eqnarray*}

We now deal with the first term of the sum on the right hand side of the last expression:
\begin{eqnarray*}
\sum_{i=1}^n\eps_i g(\nabla_{\nabla u}\nabla_{e_i} \nabla u,e_i )&=&\sum_{i=1}^n\eps_i\Big((\nabla u) g(\nabla_{e_i}\nabla u,e_i)-g(\nabla_{e_i}\nabla u,\nabla_{\nabla u}e_i)\Big)\\
&=& \left(\sum_{i=1}^n\eps_i ( (\nabla u) g( \nabla_{e_i}\nabla u,e_i) \right) - 0\\
&=& \nabla u \left(\sum_{i=1}^n\eps_i g(\nabla_{e_i}\nabla u,e_i) \right) \\
&=& \nabla u (div \nabla u)\\
&=& \nabla u (\Delta u)\\
&=& g(\nabla u , \nabla (\Delta u)).
\end{eqnarray*}

Finally, we have, recalling that $\nabla_{\nabla u} e_i (x)$ vanishes,
\begin{eqnarray*}
\sum_{i=1}^n\eps_i g(\nabla_{[\nabla u ,e_i]} \nabla u,e_i)&=& \sum_{i=1}^n\eps_i \nabla^2 u([\nabla u ,e_i], e_i)\\
&=&  \sum_{i=1}^n\eps_i \nabla^2 u( \nabla_{\nabla u} e_i - \nabla_{e_i} \nabla u, e_i)\\
&=& -  \sum_{i=1}^n\eps_i \nabla^2 u( \nabla_{e_i} \nabla u, e_i)\\
&=&- g(\nabla^2 u , \nabla^2 u).
\end{eqnarray*}

Finally, we get
$$\frac{1}{2}\Delta (g(\nabla u,\nabla u))=Ric^{\L}(\nabla u,\nabla u)+g(\nabla u,\nabla(\Delta u))+g(\nabla^2 u,\nabla^2 u),$$
the required formula.

\bigskip

\noindent

\begin{multicols}{2}
   Henri Anciaux \\
Universidade de S\~ao Paulo, IME \\
  1010 Rua do Mat\~ao,  Cidade Universit\'aria   \\
 05508-090 S\~ao Paulo, Brazil   \\
henri.anciaux@gmail.com \\

\columnbreak

Nikos Georgiou \\
Department of Mathematics and Statistics \\
University of Cyprus \\
1678 Nicosia, Cyprus \\
 georgiou.g.nicos@ucy.ac.cy

\end{multicols}

\end{document}